\documentclass[11pt, reqno]{amsart}
\usepackage{amsmath}
\usepackage{amssymb}
\newtheorem{thm*}{Theorem}
\newtheorem{thm}{Theorem}

\newtheorem{cor*}{Corollary}

\title{A characterization of the modular units}
\author{Amanda Folsom}

\begin{document}
\baselineskip = 14 pt
\parskip = 2 pt

\begin{abstract}
We provide an exact formula for the complex exponents
in the modular product expansion of the modular units, and deduce a characterization of the modular units in
terms of the growth of these exponents, answering a question of W. Kohnen.
\end{abstract}
\maketitle
\section{Introduction}
\label{intro}
Let $\varPhi(\tau)$ be the modular product defined by
\begin{eqnarray}\label{borprod} \varPhi(\tau) = \kappa q^\beta \prod_{n=1}^{\infty}
(1-q^n)^{c(n)}
\end{eqnarray}
where $q=e^{2\pi i \tau}$, $\tau \in \mathcal H$, $\beta \in
\mathbb Q$, and $c(n),\kappa \in \mathbb C$.  In
\cite{borcherds1}, Borcherds shows a duality between the exponents
$c(n)$ of the modular products $\varPhi(\tau)$  and the
coefficients $a(j)$ in the Fourier expansions $f(\tau)=\sum_{k\geq
m}a(k)q^m$ of weakly holomorphic modular forms of prescribed
weight and level.  In \cite{brunier} the authors provide an exact
formula for the exponents $c(n)$ when $\varPhi(\tau)$ is a weight
$k$ meromorphic modular form on $\Gamma(1):=\text{SL}_2(\mathbb
Z)$ whose first Fourier coefficient is $1$, in terms of the unique
modular functions $j_d(\tau)$, $d\in \mathbb N$, holomorphic on
the upper half plane $\mathcal H$ with Fourier expansion
\begin{eqnarray*} j_d(\tau)&=&q^{-d}+\sum_{m=1}^\infty
a_d(m)q^m\!\!\!\!_.
\end{eqnarray*}  The exact expansion in \cite{brunier} in this setting is given
by
\begin{eqnarray}\label{cnform}
c(n)&=&2k + n^{-1}\!\!\!\!\!\!\sum_{\tau \in \Gamma(1) \backslash
\mathcal
H^*}\!\!\!\!\!e_\tau\text{ord}_\tau(\varPhi)\sum_{d|n}\mu(n/d)j_d(\tau),
\end{eqnarray}
where the numbers $e_\tau$ are defined by $$e_\tau=\left\{%
\begin{array}{ll}
    1/2 & \hbox{$\tau = i$} \\
    1/3 & \hbox{$\tau = (1+i\sqrt{3})/2$} \\
    1 & \hbox{otherwise }, \\
\end{array}%
\right. $$ $\text{ord}_\tau(\varPhi)$ refers to the order of
$\varPhi$ at $\tau$, $\mu(n)$ is the usual
M$\ddot{\textnormal{o}}$bius function, \ \\

{\scriptsize{$^1$ Amanda Folsom, 
              Max-Planck-Institut f\"ur Mathematik, Bonn, Germany \\
	      ${\vspace{-.5in}}$
             ${\hspace{.26in}}$\texttt{alfolsom@mpim-bonn.mpg.de}}} \ \\ \ \\
	      and $^*$ denotes the
compactification of the quotient space
$\Gamma(1)\backslash\mathcal H$.
The modular units may be characterized as those meromorphic
modular functions (meromorphic modular forms of weight $0$) with divisors supported in the
cusps.  We consider the product expansions (\ref{borprod}) of the
modular units, and provide an exact formula for the exponents
$c(n)$, and from this deduce a characterization of the modular
units in terms of the growth of the exponents, answering a question of W.
Kohnen.

In what follows, for a ring $R$, we let $R^*$ denote the
multiplicative group of $R$, we let ${\mathbb Z}_\ell = \mathbb Z
/ \ell \mathbb Z$, $\mathbb Z_\ell^C = \mathbb Z_\ell \setminus
\mathbb (Z_\ell^* \cup \{0\})$, and for an integer $n$ let
$\overline{n}$ denote the equivalence class of $n$ modulo $\ell$.
We let $q_\ell:=q^{1/\ell},$ $\ell \geq 1,$ and consider modular
units of level $\ell$, that is, modular units with respect to the
principal congruence subgroups $\Gamma(\ell):=\{\gamma \in
\Gamma(1) \ | \ \gamma\equiv 1 \mod \ell \}$.

\begin{thm}\label{thm2y}
Let $u(\tau)=\varPhi(\tau/\ell)$ be a modular unit of level
$\ell=p^f$, $p$ prime, $p\neq 2,3,$ $f\in \mathbb N$. Then
\begin{eqnarray}\label{cform}c(n)&=&\frac{1}{n}\sum_{d|n}\mu(d)\sum_{k|\frac{n}{d}}t_k\Big(\frac{n}{dk}\Big)\end{eqnarray} where
\begin{eqnarray*} t_m(n) &=& \left\{%
\begin{array}{ll}
    {\displaystyle{ \ \ \ n\sum_{s\in \mathbb Z_\ell}m(\overline{n},s)e(\epsilon(n)ms/\ell)}} &
    \hbox{ \ \ \ $(n,p)=1$} \\ \ \\
    {\displaystyle{ \ \ \ n\sum_{s\in \mathbb Z_\ell^*}m(\overline{n},s)e(\epsilon(n)ms/\ell)}} &
    \hbox{ \ \ \ $p|n, \ell\nmid n$} \\ \ \\
    {\displaystyle{2n\!\!\!\!\!\sum_{s\in{\mathbb Z_\ell^*}/\{\pm 1 \}}\!\!\!\!\!\!\!m(\overline{0},s) \cos(ms/\ell)}} & \hbox{ \ \ \ $\ell |
    n$},
\end{array}%
\right.
\end{eqnarray*} $e(z)=e^{2\pi i z}$,
\begin{eqnarray}\epsilon(n) =\left\{ \begin{array}{rl} 1 & \text{if} \ n \equiv j
\!\!\!\mod \ell  \ \text{for some } j, \ 1\leq j \leq (\ell-1)/2
\\ -1 & \text{if} \ n \equiv j \!\!\!\mod \ell \ \text{for some } j,
\ (\ell-1)/2 < j \leq \ell-1,
\end{array}\right.\label{epsdef}\end{eqnarray} and
$\{m_a=m(r,s)\}_{a \in \overline{T^*_\ell}}$ is a set of integers
indexed by
\begin{eqnarray}\label{tdefn}\overline{T^*_\ell} = \{ a=(a_1,a_2)=(r/\ell,s/\ell)\in
\frac{1}{\ell}\mathbb Z^2 / \mathbb Z^2 \ | \
\textnormal{ord}(a)=\ell \}\end{eqnarray} satisfying
\begin{eqnarray}
&{\displaystyle{\sum_{a\in \overline{T^*_\ell}}m_ar^2 \equiv
\sum_{a\in \overline{T^*_\ell}}m_as^2\equiv \sum_{a\in
\overline{T^*_\ell}}m_ars \equiv 0 \!\!\!\!\mod \ell}}&
\label{quadl1}
\\ \ \nonumber \\
&\text{and}& \nonumber \\ \ \nonumber \\
&{\displaystyle{\sum_{a\in \overline{T^*_\ell}}m_a\equiv 0
\!\!\!\!\mod 12}}& \label{quadl2}
\end{eqnarray}
where $\textnormal{ord}(a)=\min\{n\in \mathbb Z^{\geq 0} \ | \
n\cdot a \in \mathbb Z^2\}$.
\end{thm}

\begin{thm}\label{thm1}  A meromorphic modular form
$u(\tau)=\varPhi(\tau/\ell)$ of weight zero on $\Gamma(\ell)$ is a
modular unit if and only if \begin{eqnarray*}c(n)\ll_u(\log \log
n)^2 \end{eqnarray*} for all $n\geq 1$, where the implied constant
depends only on $u(\tau)$.
\end{thm}

Theorem \ref{thm1} is proved more generally by W. Kohnen in \cite{kohnen} for
modular forms of weight $k$ on congruence subgroups
$\Gamma\subseteq\Gamma(1)$.  Following the statement of the
theorem in \cite{kohnen} (Theorem 1 p. 66) the author remarks 

\begin{quote}\textit{``It might be interesting to investigate if [Theorem 1 p. 66
\cite{kohnen}] could also be proved [in the case weight
$k=0$ using the theory of the modular units]."} \end{quote} 

Indeed we respond to the above remark of Kohnen and apply the
theory of the modular units to give an exact formula for the
modular exponents $c(n)$ in Theorem \ref{thm2y} (not given in
\cite{kohnen}) which allows us to prove Theorem \ref{thm1}.

\section{Modular units}
\label{sec:1}
Much of the theory of the modular units has been developed by
Kubert and Lang \cite{KL}, who provide a description of the modular units in
terms of Siegel functions. The Siegel functions are defined using
Klein forms $t_a(\tau)$, $a\in \mathbb R^2$, $\tau \in \mathcal H$
and are given by
$$\mathfrak
t_a(\tau)=e^{-\eta_a(\tau){a\cdot{(\tau,1)}}/2}\sigma_a(\tau)$$
where $\sigma$ and $\eta$ are the usual Weierstrass functions.
The Siegel functions $g_a(\tau)$ are defined by
$$g_a(\tau)=\mathfrak t_a(\tau)\Delta(\tau)^{1/12}$$ where
$\Delta(\tau)$ is the discriminant function.  The modular units of
a particular level $\ell$ form a group, and a major result of
Kubert and Lang provides a description of the modular
unit groups of level $\ell=p^f$, $p$ prime, $p\neq 2,3$, $f\in
\mathbb N$, in terms of the Siegel functions.

\begin{thm}\textnormal{(Kubert, Lang \cite{KL}.)}
For prime power $\ell = p^f$, $p\neq 2,3$ prime, $f \in \mathbb
N$, the modular units of level $\ell$ consist of products
\begin{eqnarray}\label{siegel}
\prod_{a\in \overline{T^*_\ell}}g_a^{m_a}
\end{eqnarray} of Siegel functions $g_a$, where
$\{m_a\}_{a\in \overline{T^*_\ell}}$ is a set of integers
satisfying the quadratic relations $(\ref{quadl1})$ and
$(\ref{quadl2})$.
\end{thm}   We remark that choosing different representatives in
$\overline{T^*_\ell}$ changes the Siegel function by a root of
unity, so it is understood that the theorem of Kubert and Lang is
stated modulo constants.  From the $q-$product expansion for the
function $\sigma_a(\tau)$, one may obtain the $q-$product
expansion for the Siegel functions
{\small{\begin{eqnarray}\label{siegprod}
g_a(\tau)&=&-q^{\frac{1}{2}\mathbf{B}_2(a_1)}e(a_2(a_1\!-1)
/2)\prod_{n=1}^\infty(1-q^{n-1+a_1}e(a_2))(1-q^{n-a_1}e(-a_2))
\end{eqnarray}}}where $\mathbf{B}_2(x)=x^2-x+1/6$
is the second Bernoulli polynomial.  We will use this theory to
prove Theorems \ref{thm2y} and \ref{thm1}.
\section{Proofs}
We let $\widetilde{\mathbb Z_\ell^C}$ be the set of equivalence
classes $[z]$ in $\mathbb Z_\ell^C$, defined by
$$ \widetilde{\mathbb Z_\ell^C} = \{[z] \ | \ z \in \mathbb Z_\ell^C, [z]=[w] \Leftrightarrow
 z+w\equiv 0 \!\!\!\!\!\mod \ell \}$$
and let $\overline{T_\ell^*}$ be represented by
 $$\overline{T_\ell^*} \cong
 \frac{1}{\ell}\big({\mathbb Z}_\ell^*/\{\pm 1\} \times {\mathbb Z}_\ell\big) \cup
 \frac{1}{\ell}\big(\widetilde{\mathbb Z_\ell^C} \times \mathbb Z_\ell^*\big) \cup
\frac{1}{\ell}\left(\{0\} \times {\mathbb Z_\ell^*}/\{\pm
1\}\right).
$$

Let $u(\tau)$ be a modular unit of level $\ell$, $\ell = p^f$, $p$
prime, $p\neq 2,3$, $f \in\mathbb N$.  Then there exist
$\{m_a\}_{a\in\overline{T_\ell^*}}$ satisfying (\ref{quadl1}) and
(\ref{quadl2}) such that $u(\tau)$ has an expression as given in
(\ref{siegel}).  By applying the product expansion
(\ref{siegprod}) for the Siegel functions, we see that $u(\tau)$
$$u(\tau) = \xi q^\alpha \prod_a\prod_{n=1}^\infty [(1-q^{n-1 +
a_1}e(a_2))(1-q^{n-a_1}e(-a_2))]^{m_a}$$where $\alpha \in \mathbb
Q$, $\xi \in \mathbb C$.  We compute
\begin{align}
\log (\xi q_\ell^\alpha u(\tau)^{-1}) \nonumber {\hspace{-.5in}} &
  \\
&= \sum_{a\in\overline{T_\ell^*}}\sum_{n\geq 1}\sum_{m\geq1} m_a/m
\left(q_{\ell}^{\ell m(n-1+a_1)}e(ma_2)
+ q^{\ell m(n-a_1)}e(-ma_2) \right){\hspace{2in}} \nonumber\\
&=\sum_{n\geq 1} \ \sum_{m\geq 1}
\!\!\!\!\!\!\!\!\!\!\!\!\!\!\!\!\!\sum_{ \ \ \ \ \ \ \ \ \
(r,s)\in
 {\mathbb Z}_\ell^*/\{\pm 1\} \times {\mathbb Z}_\ell} \!\!\!\!\!\!\!\!\!\!\!\!\!\!
 \!\!\!m(r,s)/m\left(q_\ell^{m(\ell(n-1)+r)}e(sm/\ell)+q_\ell^{m(\ell n-r)}e(-sm/\ell)\right)\nonumber\\
&  {\hspace{.2in}} + \sum_{n\geq 1}\sum_{m\geq
1}\!\!\!\!\!\!\!\!\!\!\!\!\!\!\!\!\!\!
\sum_{ \ \ \ \ \ \ \ \ \ (r,s)\in \widetilde{\mathbb Z_\ell^C} \times \mathbb Z_\ell^* }\!\!\!\!\!\!\!\!\!\!\!\!\!\!\!\!\!\!m(r,s)/m\left(q_\ell^{m(\ell(n-1)+r)}e(sm/\ell) + q_\ell^{m(\ell n-r)} e(-sm/\ell)\right) \nonumber\\
& {\hspace{.4in}} + \sum_{n\geq 1}\sum_{m\geq 1}\!\!\!\sum_{ \ \
s\in{\mathbb Z_\ell^*}/\{\pm
1\}}\!\!m(0,s)/m\left(q_\ell^{m(n-1)\ell}e(ms/\ell)
+ q_\ell^{m\ell n} e(-ms/\ell)\right) \nonumber \\
&= \sum_{m\geq 1}\sum_{\stackrel{n\geq 1}
{(n,p)=1}}\sum_{s\in \mathbb Z_\ell}m(\overline{n},s)/m \ e(\epsilon(n)ms/\ell)q_\ell^{mn} \nonumber\\
&{\hspace{.2in}} + \sum_{m\geq 1} \sum_{\stackrel{n\geq 1}{p|n, \ \ell\nmid n}}\sum_{s\in \mathbb Z_\ell^*}
m(\overline{n},s)/m \ e(\epsilon(n)ms/\ell)q_\ell^{mn} \nonumber\\
&{\hspace{.4in}}+ \sum_{m\geq1}\sum_{n\geq1}\sum_{s\in{\mathbb
Z_\ell^*}/\{\pm 1\}}2m(\overline{0},s)/m\cos(ms/\ell)q_\ell^{\ell
m n} \nonumber \\ & {\hspace{.6in}}+ \sum_{m\geq 1}\sum_{s\in
{\mathbb Z_\ell^*}/\{\pm 1\}}m(\overline{0},s)/m \ e(ms/\ell)
\label{lastline}
\end{align}
where $\epsilon(n)$ defined as in (\ref{epsdef}).  We will apply
the theta operator $\Theta_\ell$ with respect to the parameter
$q_\ell$, defined by
\begin{eqnarray*}\Theta_\ell (f) = q_\ell \frac{d f}{d q_\ell}. \end{eqnarray*} Using
(\ref{lastline}) we find
\begin{align}
\frac{\Theta_\ell (u)}{u} & = \alpha\ell -
\sum_{m\geq 1}\sum_{\stackrel{n\geq 1}{(n,p)=1}}\!\!\!n\!\!\sum_{s\in \mathbb Z_\ell} m(\overline{n},s) e(\epsilon(n)ms/\ell)q_\ell^{mn} \nonumber \\
&{\hspace{.5in}} - \sum_{m\geq 1} \sum_{\stackrel{n\geq 1}{p|n, \ \ell\nmid n}}\!\!\!n\!\!\sum_{s\in \mathbb Z_\ell^*} m(\overline{n},s)e(\epsilon(n)ms/\ell)q_\ell^{mn} \nonumber \\
&{\hspace{.8in}}- \sum_{m\geq1}\sum_{n\geq1}2\ell n\!\!\!\!\!\sum_{s\in{\mathbb Z_\ell^*}/\{\pm 1\}}\!\!\!\!\!m(\overline{0},s) \cos(ms/\ell)q_\ell^{\ell m n} \nonumber \\
&= \alpha\ell - \sum_{m\geq 1}\sum_{n\geq 1} t_m(n)q^{mn} \nonumber \\
&= \alpha\ell - \sum_{n\geq 1}\sum_{d|n} t_{n/d}(d)q^n
\label{theta1}
\end{align}
where \begin{eqnarray*} t_m(n) &=& \left\{%
\begin{array}{ll}
    {\displaystyle{ \ \ \ n\sum_{s\in \mathbb Z_\ell}m(\overline{n},s)e(\epsilon(n)ms/\ell)}} &
    \hbox{ \ \ \ $(n,p)=1$} \\ \ \\
    {\displaystyle{ \ \ \ n\sum_{s\in \mathbb Z_\ell^*}m(\overline{n},s)e(\epsilon(n)ms/\ell)}} &
    \hbox{ \ \ \ $p|n, \ell\nmid n$} \\ \ \\
    {\displaystyle{2n\!\!\!\!\!\sum_{s\in{\mathbb Z_\ell^*}/\{\pm 1\}}\!\!\!\!\!m(\overline{0},s) \cos(ms/\ell)}} & \hbox{ \ \ \ $\ell | n$}.
\end{array}%
\right.
\end{eqnarray*}
On the other hand, $u(\tau)$ has a modular product expansion of
the form given in (\ref{borprod}).  We compute
\begin{align*}
\log\bigg(\prod_{n\geq 1}(1-q_\ell^n)^{-c(n)}\bigg) &= \sum_{n\geq
1}\sum_{m\geq 1}c(n)/m \ q_\ell^{mn}
\end{align*}
so that
\begin{align}
\frac{\Theta_\ell (u)}{u} &= \beta -\sum_{n\geq 1}\sum_{m\geq 1}
nc(n) q_\ell^{mn} \nonumber \\ &= \beta - \sum_{n\geq 1}\sum_{d |
n} dc(d)q_\ell^n. \label{thetatwo}
\end{align}  Comparing (\ref{theta1}) and (\ref{thetatwo}), we find $\alpha\ell = \beta$ and
\begin{eqnarray*}
\sum_{d | n} t_d(n/d) = \sum_{d | n} dc(d).
\end{eqnarray*}  By M\"obius inversion, we find
\begin{align*}c(n)
=\frac{1}{n}\sum_{d|n}\mu(d)\sum_{k|\frac{n}{d}}t_k\Big(\frac{n}{dk}\Big).\end{align*}
This proves Theorem \ref{thm2y}.

To prove Theorem \ref{thm1}, suppose first
$u(\tau)=\varPhi(\tau/\ell)$ is a modular unit of level $\ell =
p^f$, $p$ prime, $p\neq 2,3$, $f\in \mathbb N$. Then by Theorem
\ref{thm2y}, the modular exponents $c(n)$ are of the form given in
(\ref{cform}).  Thus,
\begin{align*}|c(n)| &\leq \frac{1}{n}
\sum_{d|n}\sum_{k|\frac{n}{d}}\Big|t_{n/dk}(k)\Big| \\
&=
\sum_{d|n}\sum_{\stackrel{k|\frac{n}{d}}{(k,p)=1}}\Big|t_{n/dk}(k)\Big| \\
 &{\hspace{.5in}} + \sum_{d|n}\sum_{\stackrel{k|\frac{n}{d}}{p|k, \ \ell\nmid k}}\Big|t_{n/dk}(k)\Big| \\
 &{\hspace{.8in}} + \sum_{d|n}\sum_{\stackrel{k|\frac{n}{d}}{\ell|k}}\Big|t_{n/dk}(k)\Big| \\
 \\&\leq \sum_{d|n}\sum_{\stackrel{k|\frac{n}{d}}{(k,p)=1}}k\sum_{s\in \mathbb Z_\ell}|m(\overline{k},s)|\\
 &{\hspace{.5in}} + \sum_{d|n}\sum_{\stackrel{k|\frac{n}{d}}{p|k, \ \ell\nmid k}}\
  k\sum_{s\in \mathbb Z_\ell^*}|m(\overline{k},s)|
  \\
 &{\hspace{.8in}} + \sum_{d|n}\sum_{\stackrel{k|\frac{n}{d}}{\ell|k}}2k\!\!\!\!\!\sum_{s\in {\mathbb Z_\ell^*}/\{\pm 1\}}\!\!\!\!\!|m(\overline{0},s)| \\
 & \leq \ell\mathcal M_u\sum_{d|n}\sum_{\stackrel{k|\frac{n}{d}}{(k,p)=1}}k + \phi(\ell)\mathcal M_u \sum_{d|n}\sum_{\stackrel{k|\frac{n}{d}}{p|k, \ \ell\nmid k}}\
  k + \phi(\ell) \mathcal M_u/2 \sum_{d|n}\sum_{\stackrel{k|\frac{n}{d}}{\ell|k}}2k \end{align*}
  where ${\displaystyle{\mathcal M_u=\max_{a\in\overline{T_\ell^*}} \{|m_a|\}}}$, so that
  \begin{align*}
  |c(n)| &\leq \ell \mathcal M_u \sum_{d|n}\sum_{k|\frac{n}{d}} k \\
  &= \ell \mathcal M_u \sum_{d|n}\sigma_1(n/d) \\
  &\leq \ell \mathcal M_u \sum_{d|n} 2d\log\log d \\
  &\leq 2 \ell \mathcal M_u \log\log n\cdot \sigma_1(n) \\
  &\leq 4\ell \mathcal M_u (\log\log n)^2
\end{align*} hence $c(n)  \ll_u (\log \log n)^2$.
Conversely, suppose $\varPhi(\tau/\ell)$ is a weight $0$ modular
form of level $\ell=p^f$, $\ell$ prime, $\ell \neq 2,3$, $f\in
\mathbb N$, with exponents $c(n)$ satisfying $c(n) \ll_u (\log
\log n)^2$.  Then for all $\tau \in \mathcal H$, by
(\ref{thetatwo}) we see that $\Theta_\ell(\varPhi(\tau/\ell))$
converges, hence has no zeros or poles in $\mathcal H$.
\ \\ \ \\
\begin{center}\textsc{Acknowledgements}  {\vspace{-.0in}}\end{center} 

The author would like to thank \"Ozlem Imamo$\overline{\textnormal{g}}$lu for her comments on this
paper.

\bibliographystyle{amsplain}

\end{document}